\documentclass[11pt]{amsart}
\usepackage{amsmath,amssymb,enumerate}
\usepackage{epsf,epsfig,amsfonts,graphicx,color,url}

\usepackage[autostyle]{csquotes}

\usepackage[pagewise]{lineno}
\newcommand{\beq}{\begin{equation}}
\newcommand{\eeq}{\end{equation}}

\usepackage[all]{xy}
\usepackage{filecontents}

\def\R{\mathbb{R}}

\numberwithin{equation}{section}

\newtheorem{theorem}{Theorem}

\newtheorem{proposition}{Proposition}[section]

\newtheorem{definition}{Definition}[section]

\newtheorem{conj}{Conjecture}

\renewcommand{\emph}[1]{{\bfseries\itshape{#1}}}

\numberwithin{figure}{section}

\def\sR{subRiemannian }

\begin{document}

\newtheorem*{backgroundtheorem}{Background Theorem}


\title[No Periodic Normal Geodesics in $J^k(\R,\R^n)$]{No Periodic Normal Geodesics in $J^k(\R,\R^n)$}  
\author[A.\ Bravo-Doddoli]{Alejandro\ Bravo-Doddoli} 
\address{Alejandro Bravo: Dept. of Mathematics, UCSC,
1156 High Street, Santa Cruz, CA 95064}
\email{Abravodo@ucsc.edu}
\begin{abstract} 
The space of $k$-jets of $n$ real function of one real variable $x$ admits  the structure of  a Carnot group, which then has an associated Hamiltonian geodesic flow. As in any Hamiltonian flow, a natural question is the existence of periodic solutions. Does the space of $k$-jets have  periodic geodesics? This study will demonstrate the integrability of \sR geodesic flow, characterize and classify the \sR geodesics in the space of $k$-jets, and show that they are never periodic. 
\end{abstract}

\maketitle

\section{Introduction}

This paper is the generalization of \cite{ABD,RM-ABD,ABD-non}: In \cite{ABD}, the space of $k$-jets of real function of a single variable was presented as a \sR manifold, the \sR geodesic flow was defined and its integablity was verified. In \cite{RM-ABD}, the \sR geodesics were classified and some of their minimizing properties were studied. In \cite{ABD-non}, the non-existence of periodic geodesics on the space of $k$-jets of a real function of a single variable was proved.

The $k$-jets space of $n$ real functions of a single real variable, denoted here by $J^k(\R,\R^n)$ or $J^k$ for short, is a $(n(k+1)+1)$-dimensional manifold endowed with a canonical rank $n+1$ distribution, i.e., a linear sub-bundle of its tangent bundle. This distribution is globally framed by $n$ vector fields, denoted by $X_1,\cdots,X_{n+1}$ in Section \ref{sec:j-k-sr}, whose iterated Lie brackets give $J^k(\R,\R^n)$ the structure of a stratified group. Declaring $X_1,\cdots,X_{n+1}$ to be orthonormal endows $J^k(\R,\R^n)$ with the structure of a \sR manifold, which is left-invariant under the group multiplication. Like any \sR structure, the geodesics are projection of the solution to a Hamiltonian system defined on $T^*J^k$, called the geodesic flow on $J^k(\R,\R^n)$.

This paper has three main goals, the following theorem is the first.
\begin{theorem}\label{main1}
The \sR geodesic flow on $J^k(\R,\R^n)$ is integrable.
\end{theorem}
 
The bijection between geodesics on $J^k(\R,\R)$ and the pairs $(F,I)$ will be generalized, module translation $F(x) \to F(x-x_0)$, where $F(x)$ is a polynomial of degree $k$ or less and $I$ is a closed interval associated to $F(x)$,  made by Monroy-Perez and Anzaldo-Meneses \cite{Monroy1,Monroy2,Monroy3}, also described in \cite{RM-ABD} (see pg. 4). 
In the present paper it will be a bijection between the geodesic in $J^k(\R,\R^n)$ and the pairs $(F,I)$, module translation $F(x) \to F(x-x_0)$, where $F(x)=(F^1(x),\cdots,F^n(x))$ is a polynomial vector of degree $k$ or less and $I$ is a closed interval associated to $F(x)$, see Definition \ref{def:hill-interval} for more detail of $I$.

In Section \ref{sec:geo-jk}, it will be described how to build a geodesic in $J^k(\R,\R^n)$ given a pair $(F,I)$ and prove the following main result.
\begin{theorem}\label{main2}
The prescription described in Section \ref{sec:geo-jk} yields a geodesic in $J(\R,\R^n)$ parameterized by
arclength.  Conversely, any arc-length parameterized geodesic in $J^k(\R,\R^n)$ can be achieved
by this prescription applied to some polynomial vector $F(x)$ of degree $k$ or less. 
\end{theorem}

$J^k(\R,\R^n)$ comes with a projection $\Pi:J^k(\R,\R^n) \to \R^{n+1}$ onto the Euclidean plane, which projects the frame $X_1,\cdots,X_{n+1}$ onto the standard coordinate frame $\{ \partial/\partial x,\partial / \partial \theta_0^1,\cdots, \partial / \partial \theta_0^1\}$ of $\R^{n+1}$, see Section \ref{sec:j-k-sr} for the meaning of the coordinates.  

%

Using Theorem \ref{main2}, the geodesic in $J^k(\R,\R^n)$ will be classified into two main families: line-geodesics and non-line-geodesics: We say that a geodesic $\gamma(t)$  is a line-geodesics if $\gamma(t)$ corresponds to a constant polynomial vector and its projection to $\R^{n+1}$ is a line. We say that a geodesic $\gamma(t)$ is a non-line-geodesic if $\gamma(t)$ corresponds to a non-constant polynomial and its Hill interval is compact. Moreover, if $I=[x_0,x_1]$, we say that a non-line-geodesic $\gamma(t)$  is $x$-periodic (or regular), if $x_0$ and $x_1$ are regular points of $||F(x)||^2$, that is, exist $L(F,I)$ such that $x(t+L(F,I)) = x(t)$. While, $\gamma(t)$  is critical if one point or both are critical points of $||F(x)||^2$; in this case the $x$-coordinate has an asymptotic behavior to the critical point and  then the $x$-coordinate has an infinite period.

The third main result is the answer to a question by Enrico Le Donne: Does $J^k(\R,\R^n)$ have periodic geodesics?
\begin{theorem}\label{the:non-period}
$J^k(\R,\R^n)$ does not have periodic normal geodesics.
\end{theorem}

Following this classification, the only candidates to be periodic are  $x$-periodic geodesics; so the focus is on non-constant vectors correspondig to $x$-periodic geodesics.

Remark 1: Viewing $J^k(\R,\R^n)$ as a Carnot group, Theorem \ref{the:non-period} is a particular case of the conjecture made by  Enrico Le Donne.
\begin{conj}
Carnot groups do not have periodic geodesics.
\end{conj}

Remark 2: In control theory a \enquote{chained normal form} is a control system that is locally diffeomorphic to the canonical distribution for $J^k(\R,\R^n)$, see \cite{trailer-system}.

\subsection{Outline of paper}

The outline of the paper is as follows. In Section \ref{sec:j-k-sr}, the $k$-th jet space $J^k(\R,\R^n)$ is presented as a \sR manifold, as well as, the notation that will be followed throughout the work. The \sR geodesic flow is defined and the proof of Theorem \ref{main1} is given. Finally, the Carnot structure of $J^k(\R,\R^n)$ is presented. In Section \ref{sec:geo-jk}, the prescription for constructing geodesic in $J^k(\R,\R^n)$ given the pair $(F,I)$ is described, the Hamilton equation are computed and  Theorem \ref{main2} is proved.  
In Section \ref{sec:non-p-proof}, the proof of Theorem \ref{the:non-period} is given.

\subsection*{Acknowledgments}

I would like to express my gratitude to Enrico Le Donne for asking us about the existence of periodic geodesics and thus posing the problem.
I would like to thank my advisor Richard Montgomery for his  invaluable help. 
 This paper was developed with the support of the scholarship (CVU 619610) from  \enquote{Consejo de Ciencia y Tecnologia}   (CONACYT).

\section{$J^k(\R,\R^n)$ as a \sR manifold}\label{sec:j-k-sr}

The $k$-jet of a smooth function $f:\R \to \R^n$ at a point $x_0 \in \R$ is its $k$-th order Taylor expansion at $x_0$. We will this encode this $k$-jet as a $(k+2)$-tuple of real numbers as follows:
\begin{equation*}
(j^kf) = (x_0,f^k(x_0),\cdots,f^1(x_0),f(x_0)) \in \R^{n(k+1)+1}.
\end{equation*}
As $f$ varies over smooth functions and $x_0$ over $\R$, these $k$-jets sweep out the $k$-jet space. $J^k(\R,\R^n)$ is diffeomorphic to $\R^{n(k+1)+1}$ and we will use the global coordinates 
\begin{equation*}
(x,u_k,\cdots,u_1,u_0) \in \R^{n(k+1)+1}.
\end{equation*}
Where, $u_i= (u_i^1,\cdots,u_i^n)$ and, if $f = u_0$, then $u_1 = du_0/dx$, and more general, $u_{i+1} = du_{i}/dx$, $j\geq 1$. These equations are rewritten into $du_0=u_1dx$, and in general, $du_{i}=u_{i+1}dx$, we see that $J^k(\R,\R^n)$ is endowed with a natural rank $(n+1)$ distribution $D \subset TJ^k$ characterized by the $nk$ Pfaffian equations
\begin{equation*}
\begin{split}
0 &= du_0-u_1dx \\
0 &= du_1-u_2dx \\
\vdots\; & = \; \vdots \\
0 &= du_k-u_{k-1}dx. \\
\end{split}
\end{equation*}
$J^k(\R,\R^n)$ has a natural definition using the coordinates $u_i$, but they do not reflect the symmetries of the dynamics, see the proof of Theorem \ref{main1} in Section \ref{sec:geo-jk}. We will introduce the alternate coordinates $\theta_i$ for $J^k(\R,\R^n)$ describes in \cite{Monroy1,Monroy2} and also introduced in \cite{RM-ABD,ABD-non}, they are exponential coordinates of the second type, see \cite{BALOGH} Section 6.2.;
\begin{equation*}
\begin{split}
\theta_0 &= u_k \\
\theta_1 &= xu_{k}-u_{k-1} \\
\vdots & = \vdots \\
\theta_k &=\frac{x^k}{k!}u_{k} -\frac{x^k}{k!}u_{k-1}dx+\cdots+(-1)^ku_0. \\
\end{split}
\end{equation*}

$D$ is globally framed by $(n+1)$ vector fields:
\begin{equation}
X_0 = \frac{\partial}{\partial x}, \qquad X_0^j = \sum_{i=0}^k \frac{x^i}{i!}\frac{\partial}{\partial \theta^j_i}\;\;\text{for} \; 1\leq j \leq n. 
\end{equation}
A \sR structure on $J^k(\R,\R^n)$ is defined by declaring these $(n+1)$ vector fields to be orthonormal. In these coordinates the \sR metric is defined by restricting $ds^2 = dx^2 + (d\theta_0^1)^2+\cdots+(d\theta_0^n)^2$ to $D$.

During this work we will use the convention $\theta_i^j$,  where $i = 0,\cdots k$ and $j = 1,\cdots,n$ , that is, $i$ is used to denotes the vector $\theta_i$ and $j$ denote the $j$-th entry of the vector $\theta_i$

\subsection{Hamiltonian} Let $(p_x,p_{\theta_0},\cdots,p_{\theta_k},x,\theta_0,\cdots,\theta_k)$ be the traditional coordinates  for the cotangent bundle $T^*J^k$, or abbreviated as $(p,q)$. Also, let  $P_{X_0}$,$P_{X_0^1},\cdots,P_{X_0^n}:T^*J^k \to \R$ be the momentum functions of the vector fields $X_0,X_0^1,\cdots,X_0^n$, in the coordinates $(p,q)$; the momentum functions are given by 
\begin{equation}
P_{X_0} = p_x, \qquad P_{X_0^j} = \sum_{i=0}^k \frac{x^i}{i!}p_{ \theta^j_i}\;\;\text{for} \; 0\leq j \leq k. 
\end{equation}
Then the Hamiltonian governing the \sR geodesic flow on $J^k(\R,\R^n)$ is
\begin{equation}\label{eq:h-jk}
H = \frac{1}{2}(P_{X_0}^2+P_{X_0^1}^2+\cdots+P_{X_0^k}^2)
\end{equation}
(see \cite{tour},  pg 8). We will see in Section \ref{sec:geo-jk} that the condition $H = 1/2$ implies that the geodesics are parameterized by arc-length.

\subsection{Proof of Theorem \ref{main1}}

\begin{proof}
The Hamiltonian $H$ does not depend on the coordinate $\theta_i^j$ because the Hamilton equations $p_{ \theta^j_i}$ is a constant of motion. Then $\{ H,p_{ \theta^j_i}\}$ is a set of $n(k+1)+1$ constants of motion that Poisson commute and they are linearly independent. 
\end{proof}

\subsection{Carnot Group structure}
The frame $\{X_0,X_0^1,\cdots,X_0^n\}$ generates $(n(k+1)+1)$-dimensional nilpotent Lie algebra, under the iterated bracket. That is, 
\begin{equation*}
X_0^1 = [X_0,X^j_0],\;\cdots,\;X^j_k = [X_0,X_j^{k-1}], \cdots 0 = [X_0,X^j_{k}], 
\end{equation*} 
all the other Lie brackets $[X_m^\ell,X_i^{j}]$ are zero. Then the frame $\{X_0,X_i^j\}$ with $0\leq i \leq k$ and $1\leq j \leq n$ forms a $n(k+1)+1$-dimensional graded nilpotent Lie algebra:
\begin{equation*}
\mathfrak{g}_k= V_1\oplus \cdots \oplus V_{k+1}, V_1 = \{ X_0, X_0^j\}, V_i = \{ X_{i-1}^j\},\; 1\leq i \leq k ,\; 1\leq j \leq n.
\end{equation*}
Like any graded nilpotent Lie algebra, this algebra has an associated Lie group which is a Carnot group $G$ w.r.t the \sR structure. We can identify $G$ with $J^k(\R,\R^n)$, using the flows of $\{X_0,X_i^j\}$. For more detail on the jets space as a Carnot group see \cite{CarnotJets}.

\section{Geodesic in $J^k(\R,\R^n)$}\label{sec:geo-jk}

This Section describes how to build a geodesic on $J^k(\R,\R^n)$: Let us formalize the definition of the interval $I$.
\begin{definition}\label{def:hill-interval}
We say that a closed interval $I$ is a Hill interval, associated to $F(x)$, if $F^2(x) < 1$ for all $x$ in the interior of $I$ and $G^2(x) = 1$ for $x$ in the boundary of $I$. Then, $i$ is compact if and only if $F(x)$ is not a constant polynomial, if $I$ is in the form $[x_0,x_1]$, $x_0$ and $x_1$ are called endpoints of the Hill interval.
\end{definition}

Consider the  Hamiltonian system of one degree of freedom  defined on the plane phase space $(p_x,x)$ and with potential $1/2 ||F(x)||^2$, in other words, a Hamiltonian function given by
\begin{equation}\label{eq:ham-f}
H_{F}(p_x,x) = \frac{1}{2}p_x^2 + \frac{1}{2} ||F(x)||^2;
\end{equation}
then, the Hamilton equations are give by
\begin{equation}\label{eq:ham-f-1}
\dot{x} = p_x, \qquad \dot{p}_x =  (\frac{dF}{dx},F(x)),
\end{equation}
where $dF/dx$ is the derivative of the polynomial vector and $(\;,\;)$ is the Euclidean dot product on $\R^{n}$. Since the Hamiltonian is autonomous, we choose $H_{F} = 1/2$; then the dynamic takes place in the point where $||F(x)||^2 \leq 1$. If $F(x)$ is not the constant polynomial vector, and $I = [x_0,x_1]$ is the Hill interval, then $\dot{x} = 0$ if and only if $x = x_0,x_1$. Moreover, $x_0$ and $x_1$ are equilibrium points, if and only if, $x_0$ and $x_1$ are  critical points of $||F(x)||^2$, in other words, $0 = (dF/dx,F(x))$. 

Having found the solution $x(t)$, next we solve
\begin{equation}\label{eq:R-n+1}
\dot{\theta}^j_0(t) = F^j(x(t)),
\end{equation}
for $\theta_0^i$. Then, $c(t) = (x(t),\theta_0(t))$ is a curve on $\R^{n+1}$ parameterized by arc-length. Finally, we solve the horizontal lift equation associated to the curve $c(t)$
\begin{equation}\label{eq:hor-eq}
\begin{split}
\dot{\theta}^j_1 &= x(t)F^j(x(t)),\\
\dot{\theta}^j_2 &= \frac{x^2(t)}{2!}F^j(x(t)),\\
\vdots \;\; &= \;\; \vdots \\
\dot{\theta}^j_k &= \frac{x^k(t)}{k!} F^j(x(t)).\\
\end{split}
\end{equation}

\subsection{Hamilton equations}

To proof Theorem \ref{main2}, we need to write down the Hamilton equations for the geodesic flow. Since the Hamiltonian function \ref{eq:h-jk} is a left invariant function on the cotangent bundle of the Lie group $G$, the 'Lie-Poisson bracket' structure can be used  for such Hamiltonian flows to find the equations, see Appendix \cite{Arnold} or chapter 4 \cite{marsden2013introduction}. That is, if $X$ and $Y$ are left invariant vector fields then
\begin{equation}
\{ P_{X}, P_{Y} \} = - P_{[X,Y]}.
\end{equation}

In this context, the Hamilton equations are read as $\dot{f} = \{f,H\}$. With the Hamiltonian of this system, they expand to 
$$\dot{f} = \{f,P_0\} P_0 + \{f,P_{X_0^1}\} P_{X_0^1} + \cdots \{f,P_{X_0^n}\} P_{X_0^n}.$$

Using $\{P_0, P_{X_0^j}\} =- P_{X_1^j}$, we see that $P_{0}$ and $P_{X_0^j}$ evolves according to the equations
\begin{equation}\label{eq:first-layer-m}
\dot{P}_0 = -P_{X_0^1}P_{X_1^1} - \cdots - P_{X_0^2}P_{X_1^2} \qquad \dot{P}_{X_0^j} = P_0 P_{X_1^j} \;\text{for}\; 1\leq j \leq n. 
\end{equation}
For $1< i < k$, we have $\{P_{X_i^j},P_0\} = P_{X_{i+1}^{j}}$ and $\{P_{X_i^j},P_{X_\ell^m}\} = 0$, so 
\begin{equation}\label{eq:last-layer-m}
\begin{split}
\dot{P}_{X_3^j} &= P_0 P_{X_2^j} \\
\dot{P}_{X_4^j} &= P_0 P_{X_3^j} \\
\vdots \;       &= \vdots \;   \\
\dot{P}_{X_{k-1}^j} &= P_0 P_{X_k^j} \\
\dot{P}_{X_k^j} &= 0 ,\\
\end{split}
\end{equation}
for all $1 \leq j \leq n$. We also compute the Hamilton equations for the coordinates $(x,\theta_0,\cdots,\theta_n)$, 
\begin{equation}\label{eq:first-layer-c}
\dot{x} = P_{0}  \qquad \dot{\theta}_i^j = \frac{x^i}{i!} P_{X_0^j}\;\text{for}\; 0\leq i \leq k \;\text{and}\; 1\leq j \leq n. 
\end{equation}

\subsection{Proof of Theorem \ref{main2}}

\begin{proof}
Let $\gamma(t)$ be a curve corresponding to the pair $(F,I)$, that is, the coordinates $x,\theta_0^i,\theta_j^i$ are solutions to the equations (\ref{eq:ham-f}), (\ref{eq:R-n+1}) and (\ref{eq:hor-eq}), we will associate to $\gamma(t)$ some momentum functions and show that they hold equations (\ref{eq:first-layer-m}) and (\ref{eq:last-layer-m}), respectively.

Let $(p_x(t),x(t))$ be the solution to the equation (\ref{eq:ham-f-1}) with $x(t)$ laying in the $I$, comparing with the geodesic equation from (\ref{eq:first-layer-c}), we define $P_0 := p_x$. In the same way, comparing the equations (\ref{eq:R-n+1}) and (\ref{eq:hor-eq}) with the Hamilton equations (\ref{eq:first-layer-c}) and (\ref{eq:last-layer-m}) for $\theta_0^j$ and $\theta_i^j$, we define $P_{X_0^j}(t) := F^j(x(t))$ and $P_{X_i^j}(t) := \frac{d^i}{dx^i} F^j(x(t))$. Then using the change rule we have
$$ \dot{P}_{X_0^j}(t) = \frac{d}{dt} F^j(x(t)) = \frac{dF^j}{dx} \dot{x} = P_{X_1^j} P_{0}, $$
which is the equation (\ref{eq:first-layer-m}). In the same way
\begin{equation}\label{eq:last-layer-m-1}
 \dot{P}_{X_i^j}(t) = \frac{d}{dt} \frac{d^iF^j}{dx^i}(x(t)) = \frac{d^{i+1} F^j}{dx^{i+1}} \dot{x} = P_{X_1^j} P_{0}.
\end{equation}
Since $F^i(x)$ is a polynomial of degree $k$ or less, we obtain $\dot{P}_{X_k^j}(t) = 0$ for all $j = 1,\cdot, n$, and the equation (\ref{eq:last-layer-m-1}) is the same as equation (\ref{eq:last-layer-m}).

Conversely, let $\gamma(t)$ be a geodesic parameterized by arc-length with the initial condition $\gamma(0)$, that is, $\gamma(t)$ is the projection to the solution $(p(t),\gamma(t))$ of the Hamiltonian function (\ref{eq:h-jk}), we will show that the coordinates $x,\theta_0^j,\theta_i^j$ of the geodesic $\gamma(t)$ hold the equations (\ref{eq:ham-f}), (\ref{eq:R-n+1}) and (\ref{eq:hor-eq}), respectively. 

Being $\gamma(t)$ a solution to the Hamilton equations $p_{\theta^i_j}(t)$ is constant, if $a_i^j := i!p_{\theta^j_i}$ and $F^j(x) := a_0^j+a_1^jx+\cdots+a_k^jx^k$ for all $1\leq j \leq n$, then, using these expressions and $x_p = P_{X_0}$, the Hamiltonian function (\ref{eq:h-jk}) became 
\begin{equation*}
H = \frac{1}{2}(P_{X_0}^2+P_{X_0^1}^2+\cdots+P_{X_0^k}^2) = \frac{1}{2}(p_x^2 + ||F(x)||^2) = H_{F}.
\end{equation*}
Thus the $x$-coordinate of the geodesic $\gamma(t)$ is a solution to the Hamiltonian system of one degree of freedom with potential $1/2||F(x)||$, defined by equation (\ref{eq:ham-f}), where the initial condition $x(0)$ lays in a Hill interval $I$ and, so does $x(t)$. In the same way, using the solution $x(t)$ and the Hamilton equation for $\theta_0^j$, that is, $\dot{\theta}_0^j = \partial H / \partial p_{\theta_0^j} = F^j(x(t))$, thus the $\theta_0^j$-coordinate of the geodesic $\gamma(t)$ is a solution to equation (\ref{eq:R-n+1}). Finally, the Hamilton equation for $\theta_i^j$, that is, $\dot{\theta}_i^j = \partial H / \partial p_{\theta_i^j} = \frac{x^i}{i!} F^i(x(t))$ is equivalent to the horizontal equation (\ref{eq:hor-eq}). Thus, $\gamma(t)$ is a geodesic corresponding to the pair $(F,I)$.
\end{proof}

\subsection{Geodesics Classification in $J^k(\R,\R^n)$}

Using the bijection between geodesics in $J^k(\R,\R^n)$ and the pair $(F,I)$, the geodesics are classified. Let $\gamma(t)$ be a geodesic corresponding to $(F,I)$, as said before the first dichotomy is if the projected curve $\pi(\gamma(t)) = c(t)$ is a line or not.

\begin{itemize}
\item We say that $\gamma(t)$ is a line-geodesic if $F(x)$ is the constant polynomial vector, since equation (\ref{eq:R-n+1}) implies that the curve $c(t) = (x(t),\theta_0(t))$ in $\mathbb{R}^{n+1}$ is a line.

\item We say that $\gamma(t)$ is a non-line-geodesic if $F(x)$ is not the constant polynomial vector with Hill interval $I = [x_0,x_1]$, since equation (\ref{eq:ham-f-1}) implies that the $x$-dynamics takes place in $I$ and curve $c(t) = (x(t),\theta_0(t))$ in $\mathbb{R}^{n+1}$ is not a line.
\end{itemize}

Let $\gamma(t)$ be a non-line-geodesic corresponding to $(F,I)$, where $I = [x_0,x_1]$, the second dichotomy refers to the qualitative behavior of the $x(t)$ dynamic.

\begin{itemize}
\item We say that $\gamma(t)$ is $x$-periodic or regular, that is, exist $L(F,I)$ such that $x(t+L(F,I)) = x(t)$, if $x_0$ and $x_1$ are regular points of the potential $1/2||F(x)||^2$, if and only if, $x_0$ and $x_1$ are simple roots of $1-||F(x)||^2$, if and only if, $1-||F(x)||^2 = (x-x_0)(x_1-x)q(x)$, where $q(x)$ is not zero if $x$ is in $I$. 

\item We say that $\gamma(t)$ is critical, if one or both endpoints $x_0$ and $x_1$ are critical points of the potential $1/2||F(x)||^2$, if and only if, one or both endpoints $x_0$ and $x_1$ are not simple roots of $1-||F(x)||^2$. Then, by equation (\ref{eq:ham-f}), the critical points are equilibrium points of a one degree of freedom system, and the solution $x(t)$ has an asymptotic behavior to the critical points. 
\end{itemize} 

\subsubsection{Periods}

$x$-periodic geodesics have the property that the change undergone by the coordinates $\theta_j^i$ after one $x$-period $L(F,I)$ is finite and does not depend on the initial point.  This is summarized in the following proposition.

\begin{proposition}\label{prop:period}
Let $\gamma(t) = (x(t),\theta_0(t),\cdots,\theta_k(t))$ in  $J^k(\R,\R^n)$ be an $x$-periodic geodesic corresponding to the pair $(F,I)$. Then the $x$-period is 
\begin{equation}\label{eq:period}
L(F,I) = 2\int_{I} \frac{dx}{\sqrt{1-||F(x)||^2}},
\end{equation}
and is twice the time it takes for the $x$-curve to cross its Hill interval exactly once. After one period, the changes $\Delta \theta_i^j := \theta_i^j(t_0+L) - \theta_i^j(t_0)$ for $i = 0,1, \dots, k$ and $j = 1,\cdots,n$ undergone by $\theta_i^j$ are given by
\begin{equation}\label{eq:the-period}
\Delta \theta_i^j(F,I) = \frac{2}{i!} \int_{I} \frac{x^i F^j(x)dx}{\sqrt{1-||F(x)||^2}}.
\end{equation}
\end{proposition}

The proof of this Proposition is equivalent to the proofs of Proposition 4.1 from \cite{RM-ABD} (pg. 13) or Proposition 2.1 from \cite{ABD-non} (pg. 2). In \cite{RM-ABD} an argument of classical mechanics was used, see \cite{Landau} pg. 25 equation (11.5); while, in \cite{ABD-non}, a generating function to find action-angle coordinates for Hamiltonian systems was constructed, see \cite{Arnold} Section 50.  

Then a $x$-periodic geodesic $\gamma(t)$ corresponding to the pair $(F,I)$ is periodic if and only if $\Delta \theta_i^j(F,I) = 0$ for all for $i = 0,1, \dots, k$ and $j = 1,\cdots,n$.

\section{Proof of Theorem \ref{the:non-period}}\label{sec:non-p-proof}

Because that period $L(F,I)$ in equation (\ref{eq:period}) is finite, we can define an inner product in the space of polynomials of degree $k$ or less as follows
\begin{equation}
 <P_1(x),P_2(x)>_{F} := \int_{I} \frac{P_1(x)P_2(x)dx}{\sqrt{1-F^2(x)}}. 
\end{equation}
This inner product is not degenerated and will be the key to the proof of Theorem \ref{the:non-period}.

\subsection{Proof of Theorem \ref{the:non-period}}

\begin{proof}
It will be proceeded by contradiction. Let us assume $\gamma(t)$ is a periodic geodesic on $J^k(\R,\R^n)$ corresponding to the pair $(F,I)$, where $F(x)$ is not a constant polynomial vector; then $\Delta \theta_i^j(F,I) = 0$ for all $i = 0,\cdots,k$ and $j = 1,\cdots,n$. 

In the context of the space of polynomials of degree $k$ or less with inner product $< \;, \;>_{F}$, the condition $\Delta \theta_i^j(F,I) = 0$ for all $i$ and $j$ is equivalent to each $F^j(x)$ being perpendicular to $x^i$ for all $i \in 0,1,\cdots,k$ ($0 = \Delta \theta_i^j(F,I) =  <x^i,F^j(x)>_{F}$). But $\{x^i\}$ is a basis for the space of polynomials of degree $k$ or less, then each $F^j(x)$ is perpendicular to any vector, so each $F^j(x)$ is zero since the inner product is not degenerated.
This is a contradiction to the assumption that $F(x)$ is not a constant polynomial. 
\end{proof}

\end{document}